\documentclass[a4paper,11pt]{amsart}

\usepackage{amsthm, amsmath, latexsym, amsfonts}
\usepackage[all]{xy}
\usepackage{amssymb}

\theoremstyle{plain}
\newtheorem{Theorem}{Theorem}
\newtheorem{Corollary}{Corollary}
\newtheorem{Proposition}{Proposition}
\newtheorem{Lemma}{Lemma}
\newtheorem{Conjecture}{Conjecture}
\theoremstyle{definition}
\newtheorem{Definition}{Definition}

\newtheorem{Remark}{Remark}

\newcommand{\cO}{{\mathcal O}}

\newcommand{\N}{\mathbb N}
\newcommand{\Q}{\mathbb Q}

\newcommand{\R}{\mathbb R}
\newcommand{\proj}{\mathbb P}

\newcommand{\ra}{\rightarrow}

\title[Effective Adjunction]{Effective Adjunction Theory}
\author{Marco Andreatta and Claudio Fontanari}

\email{marco.andreatta@unitn.it, claudio.fontanari@unitn.it}\curraddr{
{\sc Dipartimento di Matematica \\  Universit\`a degli Studi di Trento\\
Via Sommarive 14 \\ 38123 Trento \\ Italy.}}

\thanks{We would like to thank Paolo Cascini, Roberto Pignatelli and Luis Sola-Conde for fruitful conversations.
We are grateful to J\'anos Koll\'ar for pointing out his examples and for suggesting projective varieties with canonical singularities as a good category to settle our results. 
We also thank the referees for useful comments. 
The research project was partially supported by GNSAGA of INdAM, by PRIN 2015 
"Geometria delle variet\`a algebriche", and by FIRB 2012 "Moduli spaces and Applications".}

\begin{document}

\begin{abstract}
Here we investigate the property of effectivity for adjoint divisors.
Among others, we prove the following results: 

A projective variety $X$ with at most canonical singularities is uniruled if and only if 
for each very ample Cartier divisor $H$ on $X$ we have $H^0(X, m_0K_X+H)=0$ for some $m_0=m_0(H)>0$.

Let $X$ be a projective $4$-fold, $L$ an ample divisor and $t$ an integer with $t \ge 3$. 
If $K_X+tL$ is pseudo-effective, then $H^0(X, K_X+tL) \ne 0$.
\end{abstract}

\maketitle

\section{Introduction}
Let $X$ be a normal projective variety over the complex field $\mathbb{C}$; let $K_X$ be its canonical divisor. 
We assume that $X$ has at most canonical singularities.

In the paper we fix a suitable Cartier divisor $H$ on $X$ and we discuss when the effectivity or non-effectivity of some 
adjoint divisors $aK_X + bH$ determines the geometry of $X$.

\smallskip
In the first part we consider the notion of \emph{Termination of Adjunction}. This turns out to be rather delicate, 
since in the literature there are different meanings for such a property. The following are four possibilities, where $m_0$ and $m$ are natural numbers.

\begin{itemize}

\item[(A)] For every (for some) big Cartier divisor $H$ there exists $m_0=m_0(H)>0$ such that $mK_X+H \notin \overline {Eff(X)}$ 
(i.e. it is not pseudo-effective) for $m \geq m_0$.   

\item[(B)] For every big Cartier divisor $H$ we have $H^0(X, m_0K_X+H)=0$ for some $m_0=m_0(H)>0$.
 
\item[(C)] For every very ample Cartier divisor $H$ we have $H^0(X, m_0K_X+H)=0$ for some $m_0=m_0(H)>0$.

\item[(D)] For some (for every) big Cartier divisor $H_0$ we have $H^0(X, m_0K_X+kH_0)=0$ for every $k>0$ and some $m_0=m_0(k)>0$.

\end{itemize}

It is clear that (A) $\Longrightarrow$ (B) $\Longrightarrow$ (C) $\Longrightarrow$ (D). 

\smallskip
We prove that these four definitions are equivalent and moreover that {\sl Adjunction Terminates 
in the above sense if and only if $X$ is uniruled} (see Theorem \ref{termination}, 
Corollary \ref{iff1} and Corollary \ref{iff2}).

The results follow by some characterizations of pseudo-effective Cartier divisor (see Theorem \ref{pseudo-effective}),
which are direct consequences of a fundamental result of Siu (\cite{S}).
The connection with uniruledness follows in turn from the fact that a projective variety $X$ with canonical singularities is uniruled if and only if $K_X$ is not pseudo-effective (see \cite{BDPP}, Corollary~0.3, or \cite{BCHM}, Corollary~1.3.3). 

\smallskip
A characterization of rationally connected manifolds along the same lines has been given in \cite{CDP}.

\smallskip
The examples described in \cite{Ko}, Theorem 39, show that, for varieties with singularities worst then canonical, uniruledness is not connected to Termination of Adjunction.

\medskip
We consider also the following more general definition.

\begin{itemize}

\item[(C')] Let $H$ be an effective Cartier divisor on $X$.
We say that {\sl Adjunction Terminates in the classical sense} for $H$ if 
there exists an integer $m_0 \ge 1$ such that 
$$H^0(X, H+mK_X)=0$$
for every integer $m \ge m_0$.
\end{itemize}

We conjecture that such a definition is actually equivalent to the previous ones; a partial result 
in this direction is provided by Proposition \ref{veryclassical}. 
In dimension two, Castelnuovo and Enriques indeed proved that Condition (C') implies
that $X$ is uniruled (see \cite{CE} and also \cite{Mu}). 

\bigskip
In the second part of the paper we assume that $X$ is a projective variety of dimension $n$ 
with at most terminal $\mathbb{Q}$-factorial singularities. We take a nef and and big Cartier divisor
$L$ on $X$ and we call $(X,L)$ a quasi polarized pair.

\smallskip
The following is a straightforward consequence of Theorem D in \cite{BCHM}, see Remark \ref{cascini} at the beginning of Section \ref{pabundance}.

\begin{Proposition}\label{first}
Let $(X,L)$ be a quasi polarized pair and $t >0$. 
If $K_X+tL\in \overline {Eff(X)}$, then there exists $N \in \mathbb{N}$ 
such that $H^0(X, N(K_X+tL)) \ne 0$.
\end{Proposition}

Note that for $t=0$ the statement of the Proposition would amount to Abundance Conjecture, together with MMP.

\smallskip
The next Conjecture is an effective version of the above Proposition.

\begin{Conjecture}\label{second} 
Let $(X,L)$ be a quasi polarized pair and $t > 0$. 
If $K_X+tL \in \overline {Eff(X)}$, then $H^0(X, K_X+tL) \ne 0$.
\end{Conjecture}

The case $t=1$ is a version of the so-called Ambro-Ionescu-Kawamata conjecture, 
which is true for $n \leq 3$ (see Theorem 1.5 in \cite{H}), while for $t=n-1$ 
we recover a conjecture by Beltrametti and Sommese (see \cite{BS}, Conjecture~7.2.7). 
Note that if Conjecture \ref{second} holds for $t=1$ then it holds also for every $t > 0$. 

\smallskip
In the paper we consider the following conjecture.

\begin{Conjecture}\label{third}
Let $(X,L)$ be a quasi polarized pair and $s> 0$. 
Then $H^0(X, K_X+tL)=0$ for every integer $t$ with $1 \le t \le s$ 
if and only if $K_X+sL$ is not pseudo-effective.
\end{Conjecture}

Since $L$ is big, in particular pseudo-effective, then the \emph{if} part is obvious. Note that Conjecture \ref{third} for $s=1$ implies Conjecture \ref{second}.

We prove that {\bf Conjecture \ref{third} is true for $s = n$} (see Proposition \ref{as}); we actually show that this case happens if and only if 
the pair $(X,L)$ is birationally equivalent (via a $0$-reduction, see the definition in the next section) to the pair $(\mathbb{P}^n, \mathcal{O}(1))$.

{\bf For $s=n-1$ the conjecture was essentially proved by H\"oring}, see  \cite{H}, Theorem~1.2. We prove a slightly more explicit version of 
his result (see Proposition \ref{bs}), namely, we show that this case happens if and only if the pair $(X,L)$ is birationally equivalent 
to a finite list of pairs.

\smallskip
Finally, we focus on the case $n = 4$ (see Theorem \ref{fourfold} and Proposition \ref{fourfold2}) and we generalize previous work by Fukuma 
(\cite{F2}, Theorem 3.1).

\section{Notation and preliminaries} \label{notation}

Let $X$ be a normal complex projective variety of dimension $n$. 
We adopt \cite{KollarMori} and \cite{L1} as the standard references for our set-up. 
In particular, we denote by $Div(X)$ the group of all Cartier divisors on $X$
and by $Num(X)$ the subgroup of numerically trivial divisors.
The quotient group $N^1(X) = Div(X)/Num(X)$ is the Neron-Severi group of $X$.

In the vector space $N^1(X)_{\R} := N^1(X) \otimes \R$, whose dimension is $\rho(X) := rk N^1(X)$, we consider some convex cones.
\begin{itemize}
\item[(a)] $Amp(X) \subset N^1(X)_{\R}$ the convex cone of all {\it ample} $\R$-divisor classes; it is an open convex cone.
\item[(b)] $Big(X) \subset N^1(X)_{\R}$ the convex cone of all {\it big} $\R$-divisor classes; it is an open convex cone.
\item[(e)] $Eff(X) \subset N^1(X)_{\R}$ the convex cone spanned by the classes of all effective $\R$-divisors.
\item[(n)] $Nef(X) = \overline {Amp(X)} \subset N^1(X)_{\R}$ the closed convex cone of all {\it nef} $\R$-divisor classes.
\item[(p)] $\overline {Eff(X)} = \overline {Big(X)} \subset N^1(X)_{\R}$ the closed convex cone of all {\it pseudo-effective} $\R$-divisor classes.
\end{itemize}

The above definitions actually lean on some fundamental results like the openess of the ample and big cones, the facts that $int\{\overline {Eff(X)}\} = {Big(X)}$ and 
$Nef(X) = \overline {Amp(X)} $; for more details see \cite{L1}.

Note that $Amp(X) \subset Nef(X) \cap Big(X)$ and that there are no inclusions between $Nef(X)$ and $Big(X)$.

Note also that if $ \pi : X' \ra X$ is a birational morphism and $D$ is a Cartier divisor on $X$ then $D$ is big (resp. pseudo-effective) if and only if $\pi^* D$ is big (resp. pseudo-effective).

\medskip
We consider projective varieties with singularities of special type, as in the Minimal Model Program. For reader convenience we recall their definition
(see \cite{KollarMori}, Definition 2.11 and Definition 2.12).

\begin{Definition}\label{sing}
Let $X$ be a normal projective variety. We say that $X$ has {\sl canonical }(respectively {\sl terminal}) singularities if
\begin{itemize}
\item[i)] $K_X$ is $\Q$-Cartier, and
\item[ii)] $\nu_*\cO _{\tilde X}(m K_{\tilde X})= \cO _{X}(m K_{X})$ for one (or for any) resolution of the singularities $\nu: \tilde X \to X$ 
\end{itemize}
(respectively
\begin{itemize} 
\item[ii)]  $\nu_*\cO _{\tilde X}(m K_{\tilde X} - E)= \cO _{X}(m K_{X})$ for one (or for any) resolution of the singularities $\nu: \tilde X \to X$, where $E \subset \tilde X$ is the reduced exceptional divisor).
\end{itemize}
\end{Definition}

\smallskip
In the category of projective varieties with canonical singularities the pseudo-effectivity of the canonical bundle is a birational invariant, as noticed by Mori in  \cite{M}, (11.4.1).  He actually conjectured the following beautiful result 
(\cite{M}, (11.4.2) and (11.5)), which was proved in \cite{BDPP}, Corollary~0.3 and in \cite{BCHM}, Corollary~1.3.3. 

\begin{Theorem} \label{canpseudo}
Let $X$ be a projective variety with at most canonical singularities. Then $X$ is uniruled if and only if $K_X$ is not pseudoeffective.
\end{Theorem}

\smallskip
As for the invariance of the global sections of adjoint bundles (or of pluri-canonical bundles if $L$ is trivial) we have the following. 
\begin{Lemma} \label{canonical}
Let $\pi : Y \ra X$ be a birational morphism between projective varieties with at most canonical singularities, 
let $L$ be a Cartier divisor on $X$ and let $a, b \in \N$. Then
$$H^0(X, aK_X + b L) = H^0(Y, aK_Y + b \pi^*(L)).$$
\end{Lemma}

\proof Since $Y$ and $X$ have canonical singularities we have $\pi_* a K_Y = a K_X$.  
This is straightforward from the definition of canonical singularities and by taking
a resolution of $Y$, $\nu: Y' \ra Y$,  and $\pi \circ \nu: Y' \ra X$ as a resolution of $X$.

Since $L$ is Cartier, by projection formula it follows
$$\pi_* (a K_Y + b \pi^*(L))= \pi_* (a K_Y + \pi^*(b L)) =  \pi_*(a K_Y) + b L =  a K_X + b L;$$

by taking global sections we obtain our statement.
\qed

\section{Termination of Adjunction}

Much of this section is based on the following Lemma, which was proved in the analytic setting by Siu (see \cite{S}, Proposition 1). 
For reader convenience we provide an algebraic proof relying on \cite{L2} (see also \cite{N}, Chapter V, Corollary 1.4).

\begin{Lemma}\label{siu}
Let $X$ be a smooth projective variety of dimension $n$ and let $H$ be a very ample divisor on $X$.
If $G := (n+1)H+K_X$, then for every pseudo-effective divisor $F$ on $X$ we have $H^0(X, F+G) \ne 0$.  
\end{Lemma}

\proof Since $F$ is pseudo-effective we have that $F+H$ is big, hence there exists a positive integer $m >0$
such that $m(F+H) \sim A+E$ with $A$ ample and $E$ effective (see for instance \cite{L1}, Corollary 2.2.7). 
Let $D := \frac{1}{m}E$ and $L := F+H$, so that $L-D = \frac{1}{m}A$ is big and nef; apply
\cite{L2}, Proposition~9.4.23, to get $H^0(X, K_X+L+kH + \mathcal{I}(D)) \ne 0$. 
Since the multiplier ideal $\mathcal{I}(D)$ is an ideal of $\mathcal{O}_X$, 
it follows that $H^0(X, K_X+L+kH) \ne 0$ for every $k \ge n$, i.e. $H^0(X, K_X+F+(k+1)H) \ne 0$ as soon as $k+1 \ge n+1$. 
\qed

\medskip
The following characterization of pseudo-effective divisors is probably well-known to the specialists; 
however, we did not find it explicitly in the literature.

\begin{Theorem}\label{pseudo-effective}
Let $X$ be a smooth projective variety and let $F$ be a divisor on $X$. 
The following statements, where $m$ and $N$ denote natural numbers, are equivalent:
\begin{itemize}
\item[i)] $F\in \overline {Eff(X)}$ (i.e it is pseudo-effective).
\item[ii)] There is a big divisor $G$  such that 
$H^0(X, N(mF+G)) \ne 0$ for every $m > 0$ and for some $N >0$.
\item[iii)] There is a big divisor $G$ such that $H^0(X, mF+G) \ne 0$ for all $m > 0$.
\item[iv)] There is a very ample divisor $G$ such that $H^0(X, mF+G) \ne 0$ for all $m > 0$.
\item[v)] For every big divisor $H$ we have $H^0(X, mF+kH) \ne 0$ for all $m > 0$ and all $k \geq k_0(H)$. 
\end{itemize}
\end{Theorem}

\proof
First of all note that the implications v) $\Longrightarrow$ iv), iv) $\Longrightarrow$ iii) and iii) $\Longrightarrow$ ii)
are obvious. Moreover ii) $\Longrightarrow$ i) follows from $F \equiv \lim_{m \to + \infty} \frac{mF+G}{m}$.

The difficult part is to prove i) $\Longrightarrow$ v); for this 
we use Lemma \ref{siu} together with Kodaira's Lemma (see for instance \cite{L1}, Proposition 2.2.6). 
Namely, let $G$ be the divisor of Lemma \ref{siu}; then $H^0(X, G) \ne 0$ (just take $F = \cO_X$). 
If $H$ is a big divisor on $X$, then by Kodaira's Lemma $H^0(X, kH -G) \ne 0$ for every $k \geq k_0(H)$.
Hence
\begin{eqnarray*}
\dim H^0(X, mF+kH) &=& \dim H^0(X, mF+k_0H - G + G + (k-k_0)H)\geq \\
&\geq& \dim H^0(X, mF + (k-k_0)H +G) > 0,
\end{eqnarray*}
where the last inequality follows from Lemma \ref{siu} by taking as a pseudo-effective divisor $mF + (k-k_0)H$.
\qed

\begin{Remark}
Note that i) $\Longrightarrow$ iii) is just Lemma \ref{siu}, while
i) $\Longrightarrow$ ii) follows easily from $int\{\overline {Eff(X)}\} = {Big(X)}$; this last fact was first noticed by Mori in \cite{M}, (11.3) on p. 318. 
Indeed, let $G \in Big(X)$ and $F \in \overline {Eff(X)}$; then the set $[G, F) := \{G +m F : m\in \R^+\}$ is contained in  $int\{\overline {Eff(X)}\} =Big(X)$.
\end{Remark}

\bigskip
The next Theorem proves the equivalence of the different definitions of {\sl Termination of Adjunction} stated in the Introduction. 

\begin{Theorem}\label{termination}
Let $X$ be a projective variety with at most canonical singularities. 

The following statements, where $m$ and $m_0$ denote natural numbers, are equivalent:

(i) $X$ is uniruled (i.e. $K_X$ is not pseudo-effective).

(ii) For every big Cartier divisor $H$ there exists $m_0=m_0(H)>0$ such that $mK_X+H \notin \overline {Eff(X)}$ for  $m \geq m_0$. 

(iii) For every big Cartier divisor $H$ we have $H^0(X, m_0K_X+H)=0$ for some $m_0=m_0(H)>0$. 

(iv) For every very ample Cartier divisor $H$ we have $H^0(X, m_0K_X+H)=0$ for some $m_0=m_0(H)>0$. 

(v) For some big Cartier divisor $H_0$ we have $H^0(X, m_0K_X+kH_0)=0$ for every $k>0$ and some $m_0=m_0(k)>0$.

\end{Theorem}

\proof  (i) $\Longrightarrow$ (ii) is implied by the properties of the cone described in Section \ref{notation}; indeed, it follows
by contradiction from $K_X \equiv \lim_{m \to + \infty} \frac{mK_X+H}{m}$.

(ii) $\Longrightarrow$ (iii), (iii) $\Longrightarrow$ (iv) and (iv) $\Longrightarrow$ (v) are straightforward. 

(v) $\Longrightarrow$ (i) requires a resolution of the singularities $\nu: \tilde X \to X$. Assume by contradiction that $X$ is not uniruled. 
Therefore also $\tilde X$ is not uniruled and $K_{\tilde X}$ is pseudo-effective. If $H$ is any big Cartier divisor on $X$, then 
${\tilde H}= \nu^*(H)$ is big and by \cite{L1}, Corollary 2.2.7, we have $l{\tilde H}=A+N$ with $A$ ample and $N$ effective for some $l>0$.  
It follows that $hl{\tilde H}=hA+hN$ with $hA$ very ample for some $h>0$. 
Hence, by Lemma \ref{canonical}, for every $m_0 >0$ we have 
$\dim H^0(X, m_0K_X+(n+1)hlH) = \dim H^0(\tilde X, m_0K_{\tilde X}+(n+1)hl{\tilde H})= \dim H^0(\tilde X, (m_0-1)K_{\tilde X}+(K_{\tilde X}+(n+1)hA)+(n+1)hN) \ge 
\dim H^0(\tilde X, (m_0-1)K_{\tilde X}+(K_{\tilde X}+(n+1)hA))$.
Lemma \ref{siu} says that this last term is positive, thus contradicting our assumption.
\qed

\begin{Remark}\label{history}

Note that Mori in \cite{M}, (11.4) on p. 318, suggests that in principle (i) could have been 
stronger then (iv): \emph{We say that $X$ is $\kappa$-uniruled if $K_X$ is not 
pseudo-effective. We note that $\kappa$-uniruledness is slightly stronger than saying that adjunction terminates, i.e. 
$H^0(X, mK_X+H)=0$ for each very ample divisor $H$ and some $m=m(H)>0$}. 

\end{Remark}

The following two corollaries show that the two formulations, respectively {\sl for some} and {\sl for every}, of (A) and (D) in the Introduction are equivalent.
\begin{Corollary}\label{iff1}
Let $X$ be a projective variety with at most canonical singularities. 

The following statements, where $m$ and $m_0$ denote natural numbers, are equivalent:

(i) For every big Cartier divisor $H$ there exists $m_0=m_0(H)>0$ such that $mK_X+H \notin \overline {Eff(X)}$ for  $m \geq m_0$. 

(ii) For some big Cartier divisor $H_0$ there exists $m_0=m_0(H_0)>0$ such that $mK_X+H_0 \notin \overline {Eff(X)}$ for  $m \geq m_0$. 
\end{Corollary}

\proof It is obvious that (i) implies (ii). Conversely, if (ii) holds then $K_X$ is not pseudoeffective, hence $X$ is uniruled. 
It follows from Theorem \ref{termination} that (i) holds. 
\qed

\begin{Corollary}\label{iff2}
Let $X$ be a projective variety with at most canonical singularities. 

The following statements, where $m$ and $m_0$ denote natural numbers, are equivalent:

(i) For some big Cartier divisor $H_0$ we have $H^0(X, m_0K_X+kH_0)=0$ for every $k>0$ and some $m_0=m_0(k)>0$.

(ii) For every big Cartier divisor $H$ we have $H^0(X, m_0K_X+kH)=0$ for every $k>0$ and some $m_0=m_0(k,H)>0$.
\end{Corollary}

\proof It is obvious that (ii) implies (i). Conversely, if (i) holds then by Theorem \ref{termination} $X$ is uniruled, 
i.e. $K_X$ is not pseudoeffective. Assume by contradiction that there exist a big divisor $H$ and some $k_0>0$  such that
$H^0(X, mK_X+k_0H)\ne 0$ for every $m >0$. Then $K_X = \lim_{m \to + \infty} \frac{mK_X+k_0H}{m}$ is pseudo-effective, a contradiction.
\qed

\medskip

As pointed out by the referee, since every divisor is a difference of very ample ones, (C) is actually equivalent to
the following stronger condition. 

\begin{itemize}
\item[(C*)] For every Cartier divisor $D$ we have $H^0(X, m_0K_X+D)=0$ for some $m_0=m_0(D)>0$.
\end{itemize}

\bigskip
The following is a more general definition of {\sl Termination of Adjunction}. 

\begin{Definition} \label{classicaltermination} (Condition (C'))
Let $X$ be a normal projective variety; let $H$ be an effective Cartier divisor on $X$.
We say that {\sl Adjunction Terminates in the classical sense} for $H$ if 
there exists an integer $m_0 \ge 1$ such that 
$$H^0(X, H+mK_X)=0$$
for every integer $m \ge m_0$.
\end{Definition}

We conjecture that such a definition is actually equivalent to the previous ones. The following partial result in this direction is straightforward.

\begin{Proposition}\label{veryclassical}
Let $X$ be a projective variety with canonical singularities. 
Let $H$ be any effective divisor and assume that Adjunction Terminates in the classical sense for $H$. Then $X$ has negative Kodaira dimension. 
\end{Proposition}

\proof Recall that the Kodaira dimension of a singular variety is defined to be the Kodaira dimension of any smooth model (see for instance 
\cite{L1}, Example 2.1.5). Assume by contradiction that $X$ has non-negative Kodaira dimension, i.e. $H^0(\tilde{X}, n_0 K_{\tilde{X}}) \ne 0$ 
for some integer $n_0 \ge 1$, where $\nu: \tilde{X} \to X$ is any resolution of the singularities. Since $X$ has canonical singularities, 
from Lemma \ref{canonical} it follows that $H^0(X, n_0 K_X) = H^0(\tilde{X}, n_0 K_{\tilde{X}}) \ne 0$. Hence $H^0(X, H+ n n_0 K_X) \ne 0$ 
for every integer $n \ge 1$, contradicting the assumption that $H^0(X, H+mK_X)=0$ for $m >> 0$.   
\qed

\smallskip
Together with the standard conjecture that negative Kodaira dimension implies uniruledness (see for instance \cite{M}, (11.5) on p. 319, and \cite{BDPP}, Conjecture 0.1), 
from Proposition \ref{veryclassical} it would follow that Termination of Adjunction in the classical sense implies uniruledness.  
In dimension two such an implication holds unconditionally, as it was proved by Castelnuovo and Enriques in \cite{CE} (for a modern proof we refer to \cite{Mu}).

\bigskip
We conclude this section with a characterization of uniruled varieties which
may suggest a different way to consider (effective) termination of  adjunction. It follows as a straightforward consequence 
of Lemma \ref{siu} and the main result in \cite{BDPP}. 

\begin{Proposition}
Let $X$ be a smooth projective variety of dimension $n$ and let $H$ be a very ample divisor on $X$.
If $H^0(X, mK_X+(n+1)H)=0$ for some natural number $m \ge 1$, then $X$ is uniruled.

\end{Proposition}

\proof 
Assume by contradiction that $X$ is not uniruled, so that $K_X$ is pseudo-effective by \cite{BDPP}. Lemma \ref{siu} with $F = (m-1)K_X$
gives the sought-for contradiction.   
\qed

\smallskip
Theorem 3.1 in \cite{Di} gives a statement similar to the last proposition; there the variety is singular and $H$ is just nef and big. However  
$m >1$ and $H$ has to be multiplied by a higher number, for instance $n^2$.

\section{Quasi polarized pairs}

A  {\it quasi polarized pair} is a pair  $(X,L)$ where $X$ is a projective variety
with at most $\mathbb{Q}$-factorial terminal singularities and 
$L$ is a nef and big Cartier divisor on $X$. 
If $L$ is ample we call the pair $(X,L)$ a {\it polarized pair}.

 \smallskip
In \cite{A}, Section 4, following T. Fujita's ideas as revisited by A. H\"oring in \cite{H} and using 
the MMP developed in \cite{BCHM}, we described a MMP with scaling related to divisors of type $K_X + rL$, 
for $r$ a positive rational number. 

In particular we introduced  the {\bf $0$-reduction} of a quasi polarized pair $(X,L)$ 
(see \cite{A}, Definition 4.4) as quasi polarized pair $(X',L')$ birational to $(X, L)$ 
obtained from $(X,L)$ via a Minimal Model Program with scaling:

\centerline{$(X,L) \sim (X, \Delta) :=(X_0, \Delta_0)   
\ra ----\ra (X_s, \Delta_s) \sim (X',L'),$}
which contracts or flips all extremal rays $\R^+[C]$ on $X$  such that $L\cdot  C =0$.

At every step of the MMP given above, we have a quasi polarized variety $(X_i, L_i)$ with at most terminal $\Q$-factorial singularities.

If $\pi _i:(X_{i}, \Delta_i)   \ra (X_{i+1}, \Delta_{i+1})$ is birational then $L_i = \pi_i^*(L_{i+1})$,
while if $\pi _i:(X_i, \Delta_i)   \ra (X_{i+1}, \Delta_{i+1})$ is a flip then $L_{i}$ and $ \pi_i^*(L_{i+1})$ 
are isomorphic in codimension one.

\begin{Remark}
\label{reduction} 
By using Lemma \ref{canonical} and Hartogs theorem we deduce
$$H^0(X, aK_{X} +bL) = H^0(X', aK_{X'} +bL')$$ for $a, b\in\N.$
 \end{Remark} 
 
 \smallskip
The following has been proved in \cite{A}, Theorem 5.1 and in \cite{H1}, Proposition 1.3. 

\begin{Theorem}\label{adjunction}
Let $(X,L)$ be a quasi polarized pair.
Then  $K_X + t L$ is pseudo-effective for all $t \geq n$ unless the $0$-reduction $(X',L')$ is $(\proj^n, \cO (1))$.
Actually, $K_X + (n-1)L$ is pseudo-effective unless $(X',L')$ is one of the following pairs:

\begin{itemize}
\item $(\proj ^n, \cO(1))$,
\item $(Q, \cO(1)_{|Q})$, where $Q\subset \proj ^{n+1}$ is a quadric,
\item $C_n(\proj^2, \cO(2))$, a generalized cone over $(\proj^2, \cO(2))$,
\item $X$ has the structure of a $\proj^{n-1}$-bundle over a smooth curve $C$ and $L$ restricted to any fiber is $\cO(1)$.
\end{itemize}

Moreover, except in the above cases, $K_{X'} + (n-1)L'$ is nef.
\end{Theorem}

The {\bf first-reduction} of a quasi polarized pair $(X,L)$ (see  \cite{A}, Definition 5.5) is a quasi polarized pair $(X'',L'')$ birational to $(X, L)$ obtained from 
a $0$-reduction $(X',L')$ via a morphism $\rho: X' \ra X''$ consisting of a series of divisorial contractions to smooth points,
which are weighted blow-ups of weights $(1,1,b, \ldots, b)$ with $b \geq 1$ (see \cite{AT}, Theorem 1.1). 

\begin{Remark}
\label{reduction1}
According to \cite{A}, Proposition 5.4, we have
$$H^0(K_{X} +tL) = H^0(K_{X''} +tL'')$$ for any $0\leq t \leq n-2$.
\end{Remark} 
 
\smallskip
The following has been proved in \cite{A}, Theorem 5.7.

\begin{Theorem} \label{adjunction1}
Let $(X,L)$ 
be a quasi polarized pair.

$K_{X} + (n-2)L$ is not pseudo-effective if and only if any first-reduction  $(X'', L'')$ is either 
one of the pairs listed in the statement of Theorem \ref{adjunction} or one of the following pairs: 
\begin{itemize}
\item a del Pezzo variety, that is $-K_{X''} \sim_\mathbb{Q}(n-1)L$ with $L$ ample,
\item $(\proj ^4, \cO(2))$,
\item $(\proj ^3, \cO(3))$,
\item $(Q, \cO(2)_{|Q})$, where $Q\subset \proj ^{4}$ is a quadric,
\item $X$ has the structure of a quadric fibration over a smooth curve $C$ and $L$ restricted to any fiber is $\cO(1)_{|Q}$,
\item $X$ has the structure of a $\proj^{n-2}$-bundle over a normal surface $S$ and $L$ restricted to any fiber is $\cO(1)$, 
\item $n=3$, $X$ is fibered over a smooth curve $Z$ with general fiber $\proj^{2}$ and $L$ restricted to it 
is $\cO(2)$.
\end{itemize} 

If $K_{X} + (n-2)L$ is pseudo-effective then on any first-reduction $(X'', L'')$  the divisor 
$K_{X''} + (n-2)L''$ is nef.
\end{Theorem}

\bigskip
The following definition was given by H\"oring (see (\cite{H}, Definition~1.2).

\begin{Definition}\label{scrolldef}
A quasi polarized pair $(X,L)$ is a (generalized) scroll if $X$ is smooth and there is 
a fibration $X \ra Y$ onto a projective manifold $Y$ such 
that the general fiber $F$ admits a birational morphism 
$\tau : F \ra \mathbb{P}^m$ and that $\mathcal{O}_F(L) = \tau^*\mathcal{O}_{\mathbb{P}^m}(1)$.
A quasi polarized pair $(X,L)$ is birationally a scroll if  there is a birational morphism $\nu: X' \ra X$ such that $(X', \nu^* L)$ is a (generalized) scroll.
\end{Definition}

The next is Theorem~1.4 in \cite{H}.

\begin{Theorem}
\label{hoering0}  
Let $(X,L)$ be a quasi polarized pair. If $(X,L)$ is not birationally a scroll then $\Omega_X \otimes L$ is generically nef.
\end{Theorem}

\bigskip
A key step in the proofs of Theorem \ref{bs} and of Theorem \ref{fourfold} is the following lemma 
due to H\"oring (see \cite{H}, p. 741, Step~2 in the proof of Theorem 1.2).

\begin{Lemma}
\label{hoering}
Let $(X,L)$ be a quasi polarized pair.
Assume that $K_X+(n-2)L$ is pseudo-effective and that $K_X+(n-1)L$ is  nef and big.
Then 
$$L^{n-2} [(2(K_X^2 +c_2 (X))+6nL K_X +(n+1)(3n-2)L^2] > 0.$$
\end{Lemma}

\bigskip
We consider now a quasi polarized pair $(X,L)$ and we assume moreover that $X$ is smooth. We borrow from Y. Fukuma the following set-up for the computation of the Hilbert polynomial of $K_X+tL$.

Let
\begin{eqnarray*}
F_{0}(t)&:=&\dim H^0(X, K_X +tL),\\
F_{i}(t)&:=&F_{i-1}(t+1) - F_{i -1}(t) \hbox{  for every integer  } i \hbox{  with  } 1\leq i\leq  n.
\end{eqnarray*}

\medskip
The following statement can be easily checked by reverse induction on $b \le a$. 

\begin{Lemma}\label{induction}
Fix an integer $a \ge 1$. 
If $F_{0}(t)=0$ for every integer $t$ with $1 \le t \le a$, 
then $F_{a-b}(c)=0$ for all integers $b,c$ with $1 \le c \le b \le a$.  
\end{Lemma}

\medskip
If one defines 
$$A_i(X,L) := F_{n-i}(1)$$

then it follows easily that
\begin{align}
\label{fukform}
\dim H^0(X, K_X+tL) = \sum_{j=0}^n \binom{t-1}{n-j} A_j(X,L).
\end{align}

Moreover, by taking $a = n-i+1$ and $b=c=1$ in Lemma \ref{induction}, we obtain 
the following implication.

\begin{Corollary}\label{A2}
If $H^0(X, K_X+tL)=0$ for every integer $t$ with $1 \le t \le n-i+1$, then $A_i(X,L)=0$.
\end{Corollary}

On the other hand, by Kawamata-Viehweg vanishing theorem and Serre duality, 
we have $\dim H^0(X, K_X+tL) = \chi(X,-tL)$; therefore from the Riemann-Roch theorem we obtain 
the following explicit computations (for further details, see \cite{F0}, (2.2), 
and \cite{F1}, Proposition~3.2).

\begin{Lemma}\label{fukuma}
Let $(X,L)$ be a polarized manifold of dimension $n$ and let $g(X,L)$ denote the sectional genus of $(X,L)$. 
Then we have 
\begin{eqnarray*}
A_0(X,L) &=& L^n \\
A_1(X,L) &=& g(X,L)+L^n-1 \\
24\cdot  A_2(X,L) &=&L^{n-2} [(2(K_X^2+c_2(X))+6nL K_X +(n+1)(3n-2)L^2] \\
48\cdot  A_3(X,L) &=& (n-2)(n^2-1)L^n+ n(3n-5)K_X L^{n-1} + \\ 
& & + 2(n-1) K_X^2 L^{n-2} + 2c_2(X)(K_X+(n-1)L) L^{n-3}.
\end{eqnarray*}
\end{Lemma}

\section{Polarized Abundance}
\label{pabundance}
 
The aim of this section is to argue around the Conjectures stated in the introduction. 

We start showing that Proposition \ref{first} is a direct consequence of (the more general) Theorem D in \cite{BCHM}.

\begin{Remark} \label{cascini}
Let $(X,L)$ be a quasi-polarized variety 
and let $t$ be a positive rational number. 
Then there exists an effective $\Q$-divisor $\Delta^t$ on $X$ such that
$\Delta^t \sim_{\Q} t L $ \ and
$(X, \Delta^t)$ is Kawamata log terminal.
This is well-known to the specialists, a proof can be found in \cite{A}.
If $K_X+tL\in \overline {Eff(X)}$, then $K_X+\Delta^t \in \overline {Eff(X)}$ and 
by \cite{BCHM}, Theorem D, there exists an $\R$-divisor $D \geq 0$
such that $K_X+\Delta^t \sim_{\R} D$. That is, there exists $N \in \mathbb{N}$ 
such that $H^0(X, N(K_X+tL))>0$.
\end{Remark}

\medskip
We consider Conjecture \ref{third};
for $s=n$ we recover the following easy fact. 

\begin{Proposition} \label{as}
Let $(X,L)$ be a quasi polarized pair of dimension $n$. 
We have $H^0(X, K_X+tL)=0$ for every integer $t$ with $1 \le t \le n$ 
if and only if $K_X+ nL$ is not pseudo-effective. Moreover this is the case if and only if 
the $0$-reduction $(X',L')$ of the pair $(X,L)$ is $(\mathbb{P}^n, \mathcal{O}(1))$.
\end{Proposition}

\proof By Remark \ref{reduction} we have $H^0(X, K_{X} +tL) = H^0(X', K_{X'} +tL')$ 
for any $t\geq 0$. Hence if $H^0(X, K_X+tL)=0$ for every integer $t$ with $1 \le t \le n$ 
then from Corollary \ref{A2} it follows that $A_1(X',L')=g(X',L')+L'^n-1=0$. Since we 
have $g(X',L')=0$ and $L'^n=1$ if and only if $(X',L') = (\mathbb{P}^n, \mathcal{O}(1))$, 
the claim follows from \cite{A}, Theorem~5.1~(2). 
\qed

\medskip
Next, for $s=n-1$, the following is a slightly more explicit version of \cite{H}, Theorem~1.2;
the proof is essentially the one of \cite{H}.

\begin{Theorem}\label{bs}
Let $(X,L)$ be a quasi polarized pair of dimension $n$. 
We have $H^0(X, K_X+tL)=0$ for every integer $t$ with $1 \le t \le n-1$ 
if and only if $K_X+(n-1)L$  is not pseudo-effective. 

That is, by Theorem \ref{adjunction}, if and only if 
the $0$-reduction $(X',L')$ of the pair $(X,L)$ is one of the following: 

(i) $(\proj ^n, \cO(1))$,

(ii) $(Q, \cO(1)_{|Q})$, where $Q\subset \proj ^{n+1}$ is a quadric,

(iii) $C_n(\proj^2, \cO(2))$, a generalized cone over $(\proj^2, \cO(2))$,

(iv) $X$ has the structure of a  $\proj^{n-1}$-bundle over a smooth curve $C$ 
and $L$ restricted to any fiber $F$ is $\cO(1)$.
\end{Theorem}

\proof
Let $(X',L')$ be the $0$-reduction of the pair $(X,L)$ and let  
$(\tilde X',\tilde L')$ be its desingularization 
(namely, $\nu : \tilde X' \ra X'$ and $\tilde L'= \nu ^*(L'))$.
 
By Remark \ref{reduction} and Lemma \ref{canonical} we have 
$$H^0(X, K_{X} +tL) = H^0(X', K_{X'} +tL')= H^0(\tilde X', K_{\tilde X'} +t\tilde L')$$ 
for any $t\geq 0$.

\medskip
The \emph{if} part is obvious. In order to prove the \emph{only if} part, assume that 
$H^0(X, K_{X} +tL) = H^0(X', K_{X'}$ $+tL')= H^0(\tilde X', K_{\tilde X'} +t\tilde L') =0$ 
for every integer $t$ with $1 \le t \le n-1$. Corollary \ref{A2} implies that 
\begin{equation}\label{vanish}
A_2({\tilde X'},  \tilde L')=0. 
\end{equation}

Assume by contradiction that $(X',L')$ is not one of the pairs 
in (i), (ii), (iii), (iv); then, by Theorem \ref{adjunction}, 
$K_{X'}+(n-1)L'$ is nef.  
The required contradiction is provided by \cite{H}, Theorem 1.2.
%
%
%
%
%
%
%
%
\qed 

\smallskip
The next step $s=n-2$ should work as follows.  

\begin{Conjecture}
Let $(X,L)$ be a quasi polarized manifold of dimension $n$. 
We have $H^0(X, K_X+tL)=0$ for every integer $t$ with $1 \le t \le n-2$ 
if and only if $K_X+(n-2)L$  is not pseudo-effective, that is if and only if 
the first-reduction $(X'',L'')$ is one of the pairs $(X,L)$ listed in Theorems \ref{adjunction} 
and \ref{adjunction1}.
\end{Conjecture}

Once again, the \emph{if} part is obvious.
Conversely, from Corollary \ref{A2} it follows that $A_3(X,L)=0$, but the proof of the \emph{only if} part 
seems to be elusive.

\bigskip
From now on, we focus on the case $n=4$; here formula (\ref{fukform}) reads simply as:

\begin{eqnarray}\label{fukform2}
H^0(X, K_X+tL) =
\binom{t-1}{4} A_0(X,L) + \binom{t-1}{3} A_1(X,L) + \\ 
+\binom{t-1}{2} A_2(X,L) +\binom{t-1}{1} A_3(X,L) + \binom{t-1}{0} A_4(X,L) \nonumber
\end{eqnarray}

where 
\begin{eqnarray*}
A_1(X,L) &=& g(X,L) + L^4 - 1, \\ 
A_2(X,L) &=& \dim H^0(X, K_X+3L)-2\dim H^0(X, K_X+2L)+ \\
& & + \dim H^0(X, K_X+L),\\ 
A_3(X,L) &=& \dim H^0(X, K_X+2L)-\dim H^0(X, K_X+L), \\
A_4(X,L) &=& \dim H^0(X, K_X+L).
\end{eqnarray*}

We prove the following generalization of \cite{F2}, Theorem 3.1.

\begin{Theorem}\label{fourfold}
Let $(X,L)$ be a polarized manifold of dimension $4$ and let $t$ be an integer with $t \ge 3$. 
If $K_X+tL$ is pseudo-effective, then $H^0(X, K_X+tL) \ne 0$. In particular, 
\begin{itemize}
\item $H^0(X, K_X+tL) \ne 0$ for $t\geq 5$
\item $H^0(X, K_X+4L)=0$ if and only if $(X,L)$ is $ (\mathbb{P}^4, \mathcal{O}(1))$
\item $H^0(X, K_X+3L)=0$ if and only if $(X,L)$ is either $(Q, \cO(1)_{|Q})$, where 
$Q\subset \proj ^{5}$ is a quadric, or $X$ has the structure of a $\proj^{3}$-bundle 
over a smooth curve $C$ and $L$ restricted to any fiber is $\cO(1)$.
\end{itemize}
\end{Theorem}

\proof

Since $L$ is ample $(X,L)$ is a 0-reduction, in particular by Theorem \ref{adjunction}
we can assume that $K_X+tL$ is nef for $t \geq 4$. 
We can also assume that $K_X+3L$ is nef. Indeed, if not then $(X,L)$ is one of the exceptions listed 
in the statement of Theorem \ref{adjunction}. 
If $(X,L)$ is $(\mathbb{P}^4, \mathcal{O}(1))$ or $(Q, \mathcal{O}(1))$, 
where $Q \subset \mathbb{P}^5$ is a quadric hypersurface, then Theorem \ref{fourfold} is obvious. 
The case of a generalized cone over $(\mathbb{P}^2, \mathcal{O}(2))$ does not occur 
since $X$ is smooth, while the case of a $\mathbb{P}^3$-bundle over a smooth curve
will be considered in Proposition~\ref{scroll}. 

\smallskip
Now, assume that $\Omega_X \otimes L$ is generically nef. By using the formulas in Lemma~\ref{fukuma} 
and Miyaoka inequality as stated in \cite{H}, Corollary~2.11, with $D:=4L$, we compute: 
\begin{eqnarray*}
A_2(X,L) &\ge& \frac{1}{24}\left( 2(K_X+3L)^2 L^2 + 6 (K_X+3L) L^3 + 2L^4  \right) \\
A_3(X,L) &\ge& - \frac{1}{24} (K_X+3L) L^3.
\end{eqnarray*}
Hence from (\ref{fukform2}) and the nefness of $K_X+3L$ it follows that
$$
\dim H^0(X, K_X+tL) \ge (t-1) A_3(X,L) + \frac{(t-1)(t-2)}{2} A_2(X,L) > 0
$$
for every $t \ge 3$. 

\smallskip
Finally, assume that $\Omega_X \otimes L$ is not generically nef. By Theorem~\ref{hoering0}
and Lemma \ref{canonical} we may assume that $X$ is a (generalized) scroll and the 
claim is a consequence of the following proposition. 
\qed

\begin{Proposition}\label{scroll}
Let $(X,L)$ be a generalized scroll of dimension $4$ and let $t$ be an integer such that $t \ge 3$. 
If $K_X+tL$ is nef, then $H^0(X, K_X+tL) \ne 0$. 
\end{Proposition}

\proof 
Let $X \ra Y$ be the scroll fibration and let $F$ be the generic fiber with 
a birational morphism $\tau : F \ra \mathbb{P}^m$ as in Definition \ref{scrolldef}.

\smallskip
If $X = \mathbb{P}^4$ the claim is obvious; therefore we can assume that $m \leq 3$ and that  
$A_1(X,L) = g(X,L) + L^4 - 1 > 0 $ (since we have $g(X,L)=0$ and $L^4=1$ if and only if 
$(X,L) = (\mathbb{P}^4, \mathcal{O}(1))$).
We also have that $A_0(X,L) = L^4 \geq 1$ and $A_4(X,L) = \dim H^0(X, K_X+L) \geq 0$.

\smallskip
If $m=3$, then $K_X+sL \vert_F = \tau^*\mathcal{O}_{\mathbb{P}^3}(-4+s)$, hence $H^0(X,K_X+sL)=0$ 
for $s \le 3$. 
Thus we have $A_2(X,L)=A_3(X,L)=0$ and from (\ref{fukform2}) it follows that for $t \ge 4$ 
we have
$$
\dim H^0(X, K_X+tL) \ge A_1(X,L)  > 0 
$$

\smallskip
If $m = 2$, then $K_X+sL \vert_F = \tau^*\mathcal{O}_{\mathbb{P}^2}(-3+s)$, hence $H^0(X,K_X+sL)=0$ 
for $s \le 2$. In particular, we have $A_3(X,L)=0$ and  $A_2(X,L)=\dim H^0(X, K_X+3L)$. 

For $t=3$, i.e. if we assume $K_X+3L$ is nef, by Theorem~1.2 in \cite{H} 
we must have $H^0(X,K_X+3L) \ne 0$ since $H^0(X,K_X+sL)=0$ 
for $s \le 2$.

For $t \ge 4$ we deduce from (\ref{fukform2}) that 
$$
\dim H^0(X, K_X+tL) \ge A_1(X,L)  > 0.
$$ 

\smallskip
If $m=1$, then $K_X+sL \vert_F = \tau^*\mathcal{O}_{\mathbb{P}^1}(-2+s)$, hence $H^0(X,K_X+L)=0$. 
In particular, we have $A_3(X,L) \ge 0$. 

If $H^0(X,K_X+2L)=0$, then $A_2(X,L)=\dim H^0(X, K_X+3L)$ and we 
conclude exactly as in the previous case  $m = 2$. 

If $H^0(X,K_X+2L) \ne 0$, then
$K_X+2L$ is pseudo-effective and $K_X+3L$ is pseudo-effective and big. 

Passing to the $0$-reduction we may assume that $K_X+3L$ is nef and big. Therefore Lemma~\ref{hoering} applies and by Lemma~\ref{fukuma} we get $A_2(X,L)>0$. 

Hence from 
(\ref{fukform2}) it follows that for $t \ge 3$ we have
$$
\dim H^0(X, K_X+tL) \ge A_2(X,L) > 0.
$$
\qed

\bigskip
The statement of Theorem \ref{fourfold} should hold also for $t=2$, but we have only the following partial result. 

\begin{Proposition} \label{fourfold2}
Let $(X,L)$ be a polarized manifold of dimension $4$. 
If $K_X+2L$ is pseudo-effective, then $H^0(X, K_X+2L) \ne 0$ 
unless $\Omega_X<\frac{1}{2}L>$ is not generically nef. 
\end{Proposition}

\proof 
By Theorem \ref{adjunction1} and Remark \ref{reduction1} we may assume that $K_X+2L$ is nef.

Assume that $\Omega_X<\frac{1}{2}L>$ is generically nef. By using the formula for $A_3(X,L)$ in Lemma~\ref{fukuma} 
and Miyaoka inequality, as stated in \cite{H}, Corollary~2.11, with $D:=2L$, we compute: 
$$
A_3(X,L) \ge \frac{1}{16} (K_X+2L)^2L^2 + \frac{1}{12}(K_X+2L)L^3 + \frac{1}{48} L^4.
$$
Hence from (\ref{fukform2}) it follows that
$$
\dim H^0(X, K_X+tL) \ge A_3(X,L) > 0. 
$$
\qed

\end{document}